\documentstyle[amscd,amssymb,verbatim,epsf]{amsart}

\begin{document}
\theoremstyle{plain}
\newtheorem{Thm}{Theorem}
\newtheorem{Cor}{Corollary}
\newtheorem{Ex}{Example}
\newtheorem{Con}{Conjecture}
\newtheorem{Main}{Main Theorem}
\newtheorem{Lem}{Lemma}
\newtheorem{Prop}{Proposition}

\theoremstyle{definition}
\newtheorem{Def}{Definition}
\newtheorem{Note}{Note}

\theoremstyle{remark}
\newtheorem{notation}{Notation}
\renewcommand{\thenotation}{}

\errorcontextlines=0
\numberwithin{equation}{section}
\renewcommand{\rm}{\normalshape}%

\title[Geodesic Flow]%
   {Geodesic Flow on the Normal Congruence of a Minimal Surface}
\author{Brendan Guilfoyle}
\address{Brendan Guilfoyle\\
          Department of Mathematics and Computing \\
          Institute of Technology, Tralee \\
          Clash \\
          Tralee  \\
          Co. Kerry \\
          Ireland.}
\email{brendan.guilfoyle@@ittralee.ie}
\author{Wilhelm Klingenberg}
\address{Wilhelm Klingenberg\\
 Department of Mathematical Sciences\\
 University of Durham\\
 Durham DH1 3LE\\
 United Kingdom.}
\email{wilhelm.klingenberg@@durham.ac.uk }

\date{7th March, 2006}

\begin{abstract}
We study the geodesic flow on the normal line congruence of a minimal surface in ${\Bbb{R}}^3$ induced by the 
neutral K\"ahler metric on the space of oriented lines. The metric is lorentz with isolated degenerate points and
the flow is shown to be completely integrable. In addition, we give a new holomorphic description of 
minimal surfaces in ${\Bbb{R}}^3$ and relate it to the classical Weierstrass representation. 
\end{abstract}

\maketitle

\section{Introduction}

In a recent paper \cite{gak4} a neutral K\"ahler metric was introduced on the space ${\Bbb{L}}$ of oriented affine 
lines in ${\Bbb{R}}^3$. This metric is natural in the sense that it is invariant under the action induced on 
${\Bbb{L}}$ by the Euclidean action on ${\Bbb{R}}^3$. Moreover, a surface in ${\Bbb{L}}$ is lagrangian
with respect to the associated symplectic structure iff there exist surfaces orthogonal to the associated 
2-parameter family of oriented lines (or line congruence) in ${\Bbb{R}}^3$.

In this paper we characterise the set of oriented normals to a minimal surface in ${\Bbb{R}}^3$ and study the geodesic 
flow on the line congruence induced by this neutral K\"ahler metric. Along the way, we give a new holomorphic
description of minimal surfaces in ${\Bbb{R}}^3$ and relate it to the classical Weierstrass representation. 

The induced metric on a lagrangian line congruence is either lorentz or degenerate. The null geodesics of the lorentz 
metric correspond to the principal foliation on the orthogonal surface and the degeneracy occurs precisely at 
umbilic points. We show that on the normal congruence of a minimal surface the geodesic flow is completely integrable
and find the first integrals. Recently the geodesic flow on certain non-lagrangian line congruences was investigated
\cite{gak5}. In that case, the metric was riemannian with degeneracies along a curve.

The picture that emerges is this: every minimal surface carries a completely integrable dynamical system 
\cite{guest} \cite{moser} \cite{tab}. This is generated by
geodesic motion of a lorentz metric whose null geodesics are the lines of curvatures and whose sources are the isolated
umbilic points of the minimal surface. To illustrate this we compute the geodesics explicitly for the 
case of pure harmonic minimal surfaces. These have a unique index $-N$ umbilic point (for $N>0$) and we show that the
scattering angle for non-null geodesics is $2\pi/(N+2)$.

The next section describes the normal line congruence to a 
minimal surface - all background details on the geometry of the space of oriented affine lines in ${\Bbb{R}}^3$ 
can be found in \cite{gak3} \cite{gak4} and references therein. We relate the present work to the Weierstrass 
representation in Section 3.  We then prove the result about the geodesic flow in Section 4, while we look at 
the case of pure harmonic minimal 
surfaces in the final section.

\vspace{0.2in}

\section{The Normal Line Congruence to a Minimal Surface}

Let ${\Bbb{L}}$  be the space of oriented lines in ${\Bbb{R}}^3$ which we identify with the tangent bundle to the 2-sphere
\cite{hitch}. Let $\pi:{\Bbb{L}}\rightarrow {\Bbb{P}}^1$ be the canonical bundle and (${\Bbb{J}},\Omega,{\Bbb{G}}$)
the neutral K\"ahler structure on ${\Bbb{L}}$ \cite{gak4}. 

A line congruence is a 2-parameter family of oriented lines
in ${\Bbb{R}}^3$, or equivalently, a surface $\Sigma\subset{\Bbb{L}}$. We are interested in characterising the 
line congruence formed by the oriented normal lines to a minimal surface $S$ in ${\Bbb{R}}^3$:

\begin{Thm}

A lagrangian line congruence $\Sigma\subset{\Bbb{L}}$ is orthogonal to a
minimal surface without flat points in ${\Bbb{R}}^3$ iff the congruence is the graph 
$\xi\mapsto(\xi,\eta=F(\xi,\bar{\xi}))$ of a local section of the canonical bundle with:
\begin{equation}\label{e:mineq}
\bar{\partial}\left(\frac{\partial\bar{F}}{(1+\xi\bar{\xi})^2}\right)=0,
\end{equation}
where ($\xi,\eta$) are standard holomorphic coordinates on ${\Bbb{L}}-\pi^{-1}\{\mbox{south pole}\}$
and $\partial$ represents differentiation with respect to $\xi$.
\end{Thm}
\begin{pf}
Let $S$ be a minimal surface without flat points and $\Sigma$ be its normal line congruence. Since the line congruence
is not flat, it can be given by the graph of a local section. In terms of the canonical coordinates 
($\xi,\eta=F(\xi,\bar{\xi})$) the spin coefficients of such a line congruence are \cite{gak3}:
\[
\rho=\frac{\psi}{\bar{\partial}F\partial\bar{F}-\psi\bar{\psi}}
\qquad\qquad \sigma =-\frac{\partial \bar{F}}{\bar{\partial}F\partial\bar{F}-\psi\bar{\psi}},
\]
with
\[
\psi=\partial F+r-\frac{2\bar{\xi}F}{1+\xi\bar{\xi}}.
\]
As this line congruence is orthogonal to a surface in ${\Bbb{R}}^3$, $\rho$ is real, and, as the mean curvature vanishes, 
$\rho=0$ on $S$. 

Now, the graph of a lagrangian section satisfies the following identity:
\begin{equation}\label{e:id}
(1+\xi\bar{\xi})^2\bar{\partial}\left(\frac{\sigma_0}{(1+\xi\bar{\xi})^2}\right)=-\partial\psi,
\end{equation}
where we have introduced $\sigma_0=-\partial\bar{F}$. This follows from the fact that partial derivatives commute: 
firstly the left-hand side is
\[
(1+\xi\bar{\xi})^2\bar{\partial}\left(\frac{\sigma_0}{(1+\xi\bar{\xi})^2}\right)= -\bar{\partial}\partial\bar{F}
   +\frac{2\xi\partial{\bar{F}}}{1+\xi\bar{\xi}},
\] 
while the right-hand side is
\[
-\partial\psi=-\partial\left(\bar{\partial}\bar{F}+r-\frac{2\xi\bar{F}}{1+\xi\bar{\xi}}\right)
  = -\partial\bar{\partial}\bar{F}+\frac{2\xi\partial\bar{F}}{1+\xi\bar{\xi}}.
\] 
Here we have used the lagrangian condition $\rho=\bar{\rho}$ and the equivalent local existence of a real function
$r:\Sigma\rightarrow{\Bbb{R}}$ such that
\begin{equation}\label{e:poteq}
\bar{\partial} r=\frac{2F}{(1+\xi\bar{\xi})^2}.
\end{equation}

Thus, since $\rho=0$, we have $\psi=0$ and according to the identity (\ref{e:id}), the normal congruence to a 
minimal surface must satisfy the holomorphic condition (\ref{e:mineq}).

Conversely, suppose (\ref{e:mineq}) holds for a lagrangian line congruence $\Sigma$ which is given by the graph of a 
local section. Then, by the identity (\ref{e:id}) $\psi=C$ for some real constant $C$. 
As the orthogonal surfaces move along the line congruence in ${\Bbb{R}}^3$, $\psi$ changes by
$\psi\rightarrow\psi+{\mbox{ constant}}$. Thus there exists a surface $S$ for which $\psi=0$, and therefore
$\rho=0$, i.e. there is a minimal
surface orthogonal to $\Sigma$.
\end{pf}

The previous theorem has two immediate consequences:

\begin{Cor}
The normal congruence to a minimal surface is given by a local
 section $F$ of the bundle $\pi:{\Bbb{L}}\rightarrow
S^2$ with
\[
F=\sum_{n=0}^\infty 2\lambda_n\xi^{n+3}-\bar{\lambda}_n\bar{\xi}^{n+1}
  \left({\scriptstyle{(n+2)(n+3)}}+2{\scriptstyle{(n+1)(n+3)}}\xi\bar{\xi}+{\scriptstyle{(n+1)(n+2)}}\xi^2\bar{\xi}^2
\right),
\]
for complex constants $\lambda_n$. The potential function $r:\Sigma\rightarrow{\Bbb{R}}$ satisfying (\ref{e:poteq}) is:
\[
r=-2\sum_{n=0}^\infty \frac{(3+n+(1+n)\xi\bar{\xi})(\lambda\xi^{n+2}+\bar{\lambda}\bar{\xi}^{n+2})}{1+\xi\bar{\xi}}.
\]
\end{Cor}
\begin{pf}
Since the minimal surface condition is a holomorphic condition we can expand in a power series about a point:
\[
\frac{\partial\bar{F}}{(1+\xi\bar{\xi})^2}=\sum_{n=0}^\infty \alpha_n\xi^n.
\]
This can be integrated term by term to
\[
\bar{F}=\sum_{n=0}^\infty \beta_n\bar{\xi}^n+\alpha_n\xi^{n+1}\left({\scriptstyle{\frac{1}{n+1}}}
     +{\scriptstyle{\frac{2}{n+2}}}\xi\bar{\xi}+{\scriptstyle{\frac{1}{n+3}}}\xi^2\bar{\xi}^2\right),
\]
for complex constants $\beta_n$.
Now we impose the lagrangian condition, that
\[
(1+\xi\bar{\xi})\bar{\partial}\bar{F}-2\xi\bar{F}=\sum_{n=0}^\infty \beta_n\bar{\xi}^{n-1}(n+(n-2)\xi\bar{\xi})
  -2\alpha_n\xi^{n+2}\left({\scriptstyle{\frac{1}{(n+1)(n+2)}}}+{\scriptstyle{\frac{1}{(n+2)(n+3)}}}\xi\bar{\xi}\right),
\]
is real. This implies that $\beta_0=\beta_1=\beta_2=0$ and $(n+1)(n+2)(n+3)\beta_{n+3}=-2\bar{\alpha}_n$ for $n\ge0$.
Letting $\alpha_n=-(n+1)(n+2)(n+3)\lambda_n$ gives the stated result.

Finally it is easily checked that the expressions for $r$ and $F$ satisfy (\ref{e:poteq}).
\end{pf}

On a minimal surface flat points are also umbilic points (and vice versa). Such points are now shown to be isolated:

\begin{Cor}
Umbilic points on minimal surfaces are isolated and the index of the principal foliation about an umbilic point on a
minimal surface is less than or equal to zero.
\end{Cor}
\begin{pf}
An umbilic point is a point where
$\partial\bar{F}=0$. 

Moreover, the argument of $\bar{\partial}F$ gives the principal foliation of
the surface \cite{gak3}. Given that minimality implies the holomorphic condition (\ref{e:mineq}), the zeros
of $\partial\bar{F}$ are isolated and have index greater than or equal to zero.
\end{pf}

\section{The Weierstrass Representation of a Minimal Surface}

The classical Weierstrass representation constructs a minimal surface
from a holomorphic curve in ${\Bbb{L}}$ \cite{hitch}. The minimal surface in ${\Bbb{R}}^3$ determined by a 
local holomorphic section $\nu\mapsto (\nu,w(\nu))$ of the canonical bundle is given by
\[
z=\frac{1}{2}w''-\frac{1}{2}\bar{\nu}^2\overline{w''}+\bar{\nu}\overline{w'}-\bar{w}
\]
\[
t=\frac{1}{2}\nu w''-\frac{1}{2}w'+\frac{1}{2}\bar{\nu} \overline{w''}-\frac{1}{2}\overline{w'},
\]
where a prime represents differentiation with respect to the
holomorphic parameter $\nu$ and $z=x^1+ix^2$, $t=x^3$ for Euclidean coordinates ($x^1,x^2,x^3$). The relationship
between this and our approach is as follows.

\begin{Prop}
The normal congruence of the minimal surface, in terms of the canonical
coordinates $\xi$ and $\eta$, is
\[
\xi=-\bar{\nu}    \qquad\qquad 
\eta=\frac{1}{4}(1+\xi\bar{\xi})^3\frac{\partial^2}{\partial \bar{\xi}^2}  
   \left(\frac{w}{1+\xi\bar{\xi}}\right)-\frac{1}{2}\bar{w}.
\]

\end{Prop}
\begin{pf}
We have that
\[
\frac{\partial}{\partial \nu}=\frac{1}{2}w'''
   \left(\frac{\partial}{\partial z}-\nu^2\frac{\partial}{\partial
   \bar{z}}+\nu\frac{\partial}{\partial t}\right).
\]
The unit vector in ${\Bbb{R}}^3$ which corresponds to the point
$\xi\in S^2$ is
\[
e_0=\frac{2\xi}{1+\xi\bar{\xi}}\frac{\partial}{\partial z}
   +\frac{2\bar{\xi}}{1+\xi\bar{\xi}}\frac{\partial}{\partial \bar{z}}
   +\frac{1-\xi\bar{\xi}}{1+\xi\bar{\xi}}\frac{\partial}{\partial t}.
\]
The normal direction is given by the vanishing of the inner product of
the preceding 2 vectors, which is easily seen to imply (for $w'''\neq 0$)
$\xi=-\bar{\nu}$. At $w'''=0$ there is an umbilic point.
 The remainder of the proposition follows from the incidence relation \cite{gak3}:
\[
\eta=\frac{1}{2}\left(z-2t\xi-\bar{z}\xi^2\right).
\] 
\end{pf}

The holomorphic functions of our method and that of the Weierstrass
representation are related by
\[
\frac{1}{(1+\xi\bar{\xi})^2}\frac{\partial F}{\partial \bar{\xi}}  
    =\frac{1}{4}\frac{\partial^3 w}{\partial \bar{\xi}^3}.
\]  

\vspace{0.2in}

\section{The Geodesic Flow}

We now look at the metric on lagrangian sections:

\begin{Prop}\label{p:fint}
The metric induced by the neutral K\"ahler metric on the graph of a lagrangian section $\eta=F(\xi,\bar{\xi})$
is:
\[
ds^2=\frac{2i}{(1+\xi\bar{\xi})^2}\left(\sigma_0d\xi\otimes d\xi-\bar{\sigma}_0d\bar{\xi}\otimes d\bar{\xi}\right),
\]
where $\sigma_0=-\partial \bar{F}$. Thus, for $|\sigma_0|\neq0$ the metric is lorentz and for
$|\sigma_0|=0$ the metric is degenerate. 
\end{Prop}
\begin{pf}
The neutral K\"ahler metric has local expression \cite{gak4}:
\begin{equation}\label{e:metric}
{\Bbb{G}}=\frac{2i}{(1+\xi\bar{\xi})^2}\left(
  d\eta \otimes d\bar{\xi}-d\bar{\eta}\otimes d\xi
   +\frac{2(\xi\bar{\eta}-\bar{\xi}\eta)}{1+\xi\bar{\xi}}d\xi\otimes d\bar{\xi}
\right).
\end{equation}
We pull the metric back to the section:
\[
{\Bbb{G}}|_\Sigma=\frac{2i}{(1+\xi\bar{\xi})^2}\left[
  \bar{\partial}F\;d\bar{\xi} \otimes d\bar{\xi}-\partial\bar{F}\;d\xi\otimes d\xi
   +\left(\partial F-\bar{\partial}\bar{F}+\frac{2(\xi\bar{\eta}-\bar{\xi}\eta)}{1+\xi\bar{\xi}}\right)d\xi\otimes d\bar{\xi}\right].
\]
Now the lagrangian condition says precisely that the coefficient of the $d\xi\otimes d\bar{\xi}$ term vanishes, and the 
result follows.
\end{pf}

We turn now to the geodesic flow. Since the metric above is flat on the normal congruence of a minimal surface, this
flow is completely integrable:

\begin{Prop}\label{p:main1}
Consider the normal congruence to a minimal surface $\Sigma\subset{\Bbb{L}}$ given by $(\xi,\eta=F(\xi,\bar{\xi}))$.
 The geodesic flow on $\Sigma$ is completely integrable with first integrals
\[
{\mbox {I}}_1=\frac{2i}{(1+\xi\bar{\xi})^2}\left(\sigma_0\dot{\xi}^2-\bar{\sigma}_0\dot{\bar{\xi}}^2\right)
\qquad\qquad
{\mbox {I}}_2=\frac{\sigma_0^{\scriptstyle{\frac{1}{2}}}\dot{\xi}+\bar{\sigma}_0^{\scriptstyle{\frac{1}{2}}}\dot{\bar{\xi}}}{1+\xi\bar{\xi}}.
\]
\end{Prop}
\begin{pf}
Consider the affinely parameterised geodesic $t\mapsto(\xi(t),\eta=F(\xi(t),\bar{\xi}(t)))$ on $\Sigma$ with tangent vector
\[
{\mbox T}=\dot{\xi}\frac{\partial}{\partial\xi}+\dot{\bar{\xi}}\frac{\partial}{\partial\bar{\xi}}.
\]
The geodesic equation 
${\mbox T}^j\nabla_j{\mbox T}^k=0$, projected onto the $\xi$ coordinate is
\[
\ddot{\xi}+\Gamma_{\xi\xi}^\xi\dot{\xi}^2 +2\Gamma_{\bar{\xi}\xi}^\xi\dot{\xi}\dot{\bar{\xi}}  
                       + \Gamma_{\bar{\xi}\bar{\xi}}^\xi\dot{\bar{\xi}}^2=0.
\] 
For the induced metric (as given in Proposition \ref{p:fint}) a straight-forward calculation yields the 
Christoffel symbols:
\[
\Gamma_{\xi\xi}^\xi=\frac{1}{2\sigma_0}\left(\partial\sigma_0-\frac{2\sigma_0\bar{\xi}}{1+\xi\bar{\xi}} \right)
\qquad\qquad
\Gamma_{\xi\xi}^{\bar{\xi}}=\frac{1}{2\bar{\sigma}_0}\left(\bar{\partial}\sigma_0-\frac{2\sigma_0\xi}{1+\xi\bar{\xi}} \right)
\]
\[ 
\Gamma_{\xi\bar{\xi}}^\xi=\frac{1}{2\sigma_0}\left(\bar{\partial}\sigma_0-\frac{2\sigma_0\xi}{1+\xi\bar{\xi}}\right).
\]
For the normal congruence of a minimal surface the holomorphic condition (\ref{e:mineq}) implies that 
$\Gamma_{\xi\xi}^{\bar{\xi}}=0$ and $\Gamma_{\xi\bar{\xi}}^\xi=0$.
Thus the geodesic equation reduces to
\[
\ddot{\xi}=-\frac{1}{2}\partial\left[\ln\left(\frac{\sigma_0}{(1+\xi\bar{\xi})^2}\right)\right]\dot{\xi}^2.
\]

The fact that ${\mbox I}_1$ is constant along a geodesic comes from the fact that the geodesic flow preserves 
the length of the 
tangent vector ${\mbox T}^j$. On the other hand, differentiating ${\mbox I}_2$ with respect to $t$:
\begin{align}
\dot{I}_2&=\frac{1}{2}\left(\frac{\sigma_0}{(1+\xi\bar{\xi})^2}\right)^{-\scriptstyle{\frac{1}{2}}}
  \left[\partial\left(\frac{\sigma_0}{(1+\xi\bar{\xi})^2}\right)\dot{\xi}^2
          +\bar{\partial}\left(\frac{\sigma_0}{(1+\xi\bar{\xi})^2}\right)\dot{\xi}\dot{\bar{\xi}}\right]\nonumber\\
&\qquad
   +\frac{1}{2}\left(\frac{\bar{\sigma}_0}{(1+\xi\bar{\xi})^2}\right)^{-\scriptstyle{\frac{1}{2}}}
  \left[\partial\left(\frac{\bar{\sigma}_0}{(1+\xi\bar{\xi})^2}\right)\dot{\xi}\dot{\bar{\xi}}
          +\bar{\partial}\left(\frac{\bar{\sigma}_0}{(1+\xi\bar{\xi})^2}\right)\dot{\bar{\xi}}^2\right]
   \nonumber\\
&\qquad\qquad+\frac{\sigma_0^{\scriptstyle{\frac{1}{2}}}}{1+\xi\bar{\xi}}\ddot{\xi}
        +\frac{\bar{\sigma}_0^{\scriptstyle{\frac{1}{2}}}}{1+\xi\bar{\xi}}\ddot{\bar{\xi}}\nonumber\\
&=\frac{1}{2}\left(\frac{\sigma_0}{(1+\xi\bar{\xi})^2}\right)^{-\scriptstyle{\frac{1}{2}}}
    \partial\left(\frac{\sigma_0}{(1+\xi\bar{\xi})^2}\right)\dot{\xi}^2
          +\frac{1}{2}\left(\frac{\bar{\sigma}_0}{(1+\xi\bar{\xi})^2}\right)^{-\scriptstyle{\frac{1}{2}}}
    \bar{\partial}\left(\frac{\bar{\sigma}_0}{(1+\xi\bar{\xi})^2}\right)\dot{\bar{\xi}}^2\nonumber\\
&\qquad\qquad
  -\frac{1}{2}\frac{\sigma_0^{\scriptstyle{\frac{1}{2}}}}{1+\xi\bar{\xi}}
     \partial\left[\ln\left(\frac{\sigma_0}{(1+\xi\bar{\xi})^2}\right)\right]\dot{\xi}^2
    -\frac{1}{2}\frac{\bar{\sigma}_0^{\scriptstyle{\frac{1}{2}}}}{1+\xi\bar{\xi}}
     \bar{\partial}\left[\ln\left(\frac{\bar{\sigma}_0}{(1+\xi\bar{\xi})^2}\right)\right]\dot{\bar{\xi}}^2
\nonumber\\
&=0, \nonumber
\end{align}
as claimed.
\end{pf}

\vspace{0.2in}

\section{Examples: The Pure Harmonics}

We now consider the geodesic flow for the pure harmonics, that is, the minimal surfaces with
\[
\frac{\partial\bar{F}}{(1+\xi\bar{\xi})^2}=\alpha_N\xi^N,
\]
for some $N\in{\Bbb{N}}$. These have isolated umbilic points of index $-N<0$ at $\xi=0$, which is
also an $N+1\;-$fold branch point. By a rotation 
we can make $\alpha_N$ real and rescaling the first integrals we will set it to 1. 

By Proposition \ref{p:main1} above, the first integrals are
\[
{\mbox {I}}_1=2i\left(\xi^N\dot{\xi}^2-\bar{\xi}^N\dot{\bar{\xi}}^2\right)
\qquad\qquad
{\mbox {I}}_2=\xi^{N/2}\dot{\xi}+\bar{\xi}^{N/2}\dot{\bar{\xi}}.
\]
These can be integrated to:
\[
\frac{4i{\mbox{I}}_2}{N+2}\left(\xi^{\scriptstyle{\frac{N+2}{2}}}-\bar{\xi}^{\scriptstyle{\frac{N+2}{2}}}\right)
  ={\mbox{I}}_1t+c_1
\qquad\qquad
\frac{2}{N+2}\left(\xi^{\scriptstyle{\frac{N+2}{2}}}+\bar{\xi}^{\scriptstyle{\frac{N+2}{2}}}\right)
  ={\mbox{I}}_2t+c_2,
\]
for real constants of integration $c_1$ and $c_2$.

For null geodesics ${\mbox{I}}_1=0$, and if we let $\xi=Re^{i\theta}$ we get two sets of null geodesics (future- and 
past-directed) which are given implicitly by:
\[
R^{\scriptstyle{\frac{N+2}{2}}}\sin\left(\frac{N+2}{2}\right)\theta=c_1
\qquad\qquad
R^{\scriptstyle{\frac{N+2}{2}}}\cos\left(\frac{N+2}{2}\right)\theta=c_2.
\]
For $N=0$, these form a rectangular grid, while for $N>0$ they form the standard index $-N$ foliation
about the origin. The diagram below shows the $N=1$ minimal surface, and the 
foliation of null geodesics about the index -1 umbilic.

\vspace{0.1in}
\setlength{\epsfxsize}{4.5in}
\begin{center}
   \mbox{\epsfbox{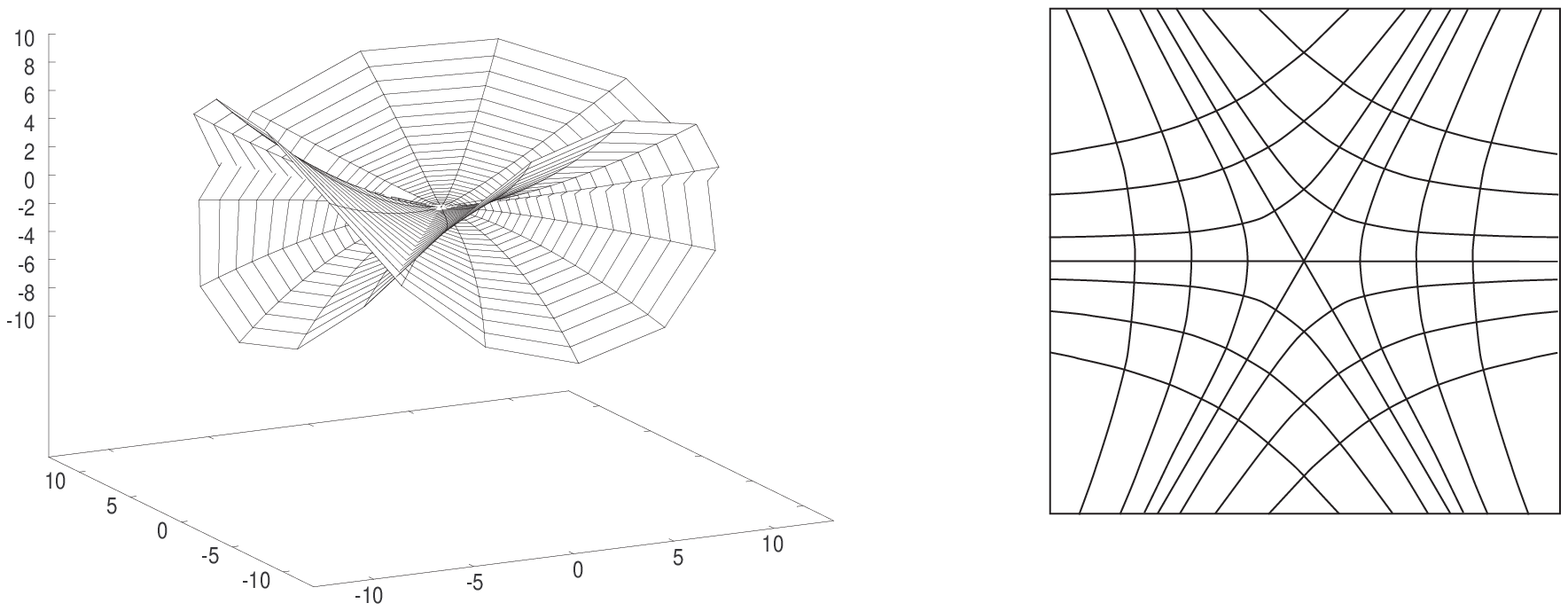}}
\end{center}
\vspace{0.1in}

For non-null geodesics ${\mbox{I}}_1\neq0$ and ${\mbox{I}}_2\neq0$. Then the geodesics can be written
parametrically:
\[
R^{\scriptstyle{\frac{N+2}{2}}}\sin\left({\scriptstyle{\frac{N+2}{2}}}\right)\theta
             =-\frac{N+2}{8{\mbox{I}}_2}\left({\mbox{I}}_1t+c_1\right)
\qquad
R^{\scriptstyle{\frac{N+2}{2}}}\cos\left({\scriptstyle{\frac{N+2}{2}}}\right)\theta
             =\frac{N+2}{4}\left({\mbox{I}}_2t+c_2\right).
\]
The umbilic acts as a source of repulsion and the scattering angle can be found by noting that
\[
\tan\left({\scriptstyle{\frac{N+2}{2}}}\right)\theta=-\frac{1}{2{\mbox{I}}_2}\frac{{\mbox{I}}_1t+c_1}{{\mbox{I}}_2t+c_2}.
\]
Thus, as $t\rightarrow\pm\infty$ we have $\tan\left({\scriptstyle{\frac{N+2}{2}}}\right)\theta\rightarrow-\frac{{\mbox{I}}_1}{2{\mbox{I}}_2^2}$. We deduce then that the scattering angle is $\frac{2\pi}{N+2}$. The diagram below illustrates the
scattering angle for a non-null geodesic about the $N=1$ umbilic.

\vspace{0.1in}
\setlength{\epsfxsize}{3.0in}
\begin{center}
   \mbox{\epsfbox{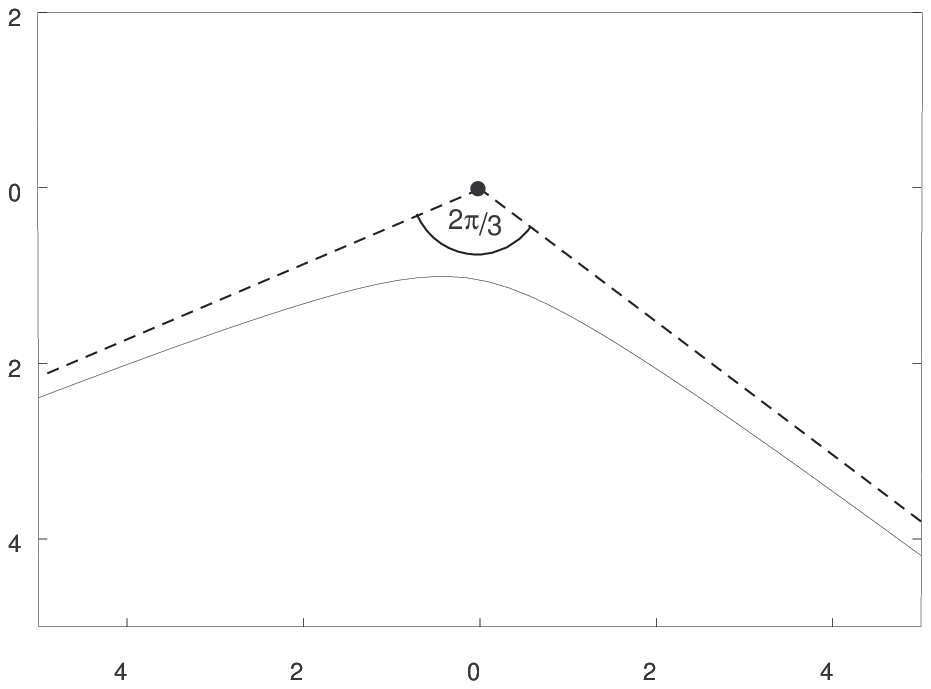}}
\end{center}

\end{document}